\input amstex
\documentstyle{amsppt}

\loadeufb
\loadeusb
\loadeufm
\loadeurb
\loadeusm

\magnification =\magstephalf
\refstyle{A}
\NoRunningHeads

\topmatter
\title  Bergman type metrics in tower of coverings \endtitle
\author Robert  Treger \endauthor
\address Princeton, NJ 08540  \endaddress
\email roberttreger117{\@}gmail.com \endemail
\endtopmatter

\document 

The aim of this note is to prove a version of a statement
attributed  to Kazhdan by Yau  \cite{Y}. For a survey of known results
and historical remarks, see the recent article by Ohsawa \cite{O, Section 5} as well as articles by Kazhdan (\cite{Ka1}, \cite{Ka2}) and McMullen \cite{M, Appendix}. As a corollary, we obtain a more transparent form of the uniformization theorem in \cite{T1}.

\head
1.  Preliminaries
\endhead

\subhead {1.1}\endsubhead Throughout the note, $X\subset \bold P^{\ell}$ denotes a nonsingular $n$-dimensional connected projective variety with large (see \cite{Kol}) and residually finite fundamental group, and without general elliptic curvilinear sections.
Let $U$ denote its universal covering. 

 We assume the canonical bundle $\Cal K_X$ on $X$ is ample. Set $\Cal K^s_X :=\Cal K^{\otimes s}_X$.  We consider a
tower of Galois coverings with each $Gal (X_i / X)$ a finite group:
$$
\qquad\qquad X=X_0 \leftarrow X_1 \leftarrow X_2 \leftarrow \cdots \leftarrow U, \quad  \bigcap_i Gal(U/X_i) = \{1\} \quad (0\leq i<\infty).
$$

It is known that for all sufficiently large $s$, the bundle $\Cal K^s_X$ and its
pullbacks on each $X_i$, denoted by $\Cal K^s_{X_i}$, are very ample (see
\cite{Kol, 16.5}). In the sequel, we fix the smallest integer $m$ with this property and assume $X\subset \bold P^{\ell}$ is given by $\Cal K_X^m$. The pullback of  $\Cal K^s_X$   on $U$ will be denoted by $\Cal K^s_U$. 

In \cite{T2, (1.1)}, we have defined a Kahler metric (called {\it original}\/ metric below) on $U$ depending on the embedding of $X \subset \bold P^{\ell}$. This metric defines the volume form $dv_U=\rho d\mu$, where $d\mu$ is the Euclidean volume form. The metric and the volume form are invariant under the action of $\pi_1(X)$ so we get the volume form $dv_U$ on $X_i$.

The bundle $\Cal K^m_U$ is equipped with the Hermitian metric $h_{\Cal K^m_U}:=h^m_{\Cal K_U}$, where $h_{\Cal K_U}$ is the standard  Hermitian metric  on $\Cal K_U$ \cite {Kol, 5.12, 5.13, 7.1}.
Similarly, each bundle $\Cal K^m_{X_i}$ is equipped with the Hermitian metric $h_{\Cal K^m_{X_i}}:=h^m_{\Cal K_{X_i}}$.

The Calabi diastasis was introduced in \cite{C} (see also \cite{U, Appendix}). For properties of infinite-dimensional projective spaces, see \cite{C, Chap.\,4} and \cite{Kob, Sect.\,7}.

\subhead 1.2.  Definition
\endsubhead
We  will define  the $q$-{\it Bergman metric}\/ on $U$, where $q\geq m$ is an integer. 
Let $H$ be  the Hilbert space of weight $q$ square-integrable with respect to $dv_U$ holomorphic differential forms $\omega$ on $U$, where by square-integrable we mean
$$
 \int_U\! h_{\Cal K^q_U}\! (\omega,\omega) \rho^{-q} dv_U < \infty. 
$$

We assume $H \not = 0$.   This Hilbert space has a {\it reproducing kernel} as, e.g., in \cite{FK, pp.\, 8-10 or pp.\, 187-189}. Let $\bold P(H^*)$ denote the corresponding projective space with its Fubini-Study metric. 
If the natural evaluation map
$$
e_U: U \longrightarrow \bold P(H^*)
$$
is an {\it embedding}\/ then the metric on $U$, induced from $\bold P(H^*)$,  is called the $q$-{\it Bergman metric}.

Similarly, one defines the Euclidean space $V_i$,  $X_i \hookrightarrow  \bold P(V_i^*)$, and the $q$-Bergman metric on $X_i$. 
We consider every $V_i$ with the {\it normalized}\/ inner product in order to get the inclusion of Euclidean spaces $V_i \hookrightarrow V_{i+1}$ induced by pullbacks of forms.

\head  2. Proposition and Corollary \endhead

\proclaim{Proposition}  With the above notation and assumptions, we assume, in addition, that $U$ has the $q$-Bergman metric for an integer $q\geq m$. Then the
$q$-Bergman metric on $U$ equals the limit of pullbacks of the $q$-Bergman metrics from $X_i$'s.
\endproclaim

\demo{Proof} (2.1) Let $b_{U,q}$ denote the
$q$-Bergman metric on $U$. Set $\tilde b_{U,q}:= \lim \sup  b_i $, where $b_i $ is the pullback on $U$ of the $q$-Bergman
metric on $X_i$.  

 By an argument similar to \cite {Ka2, Lemma}, now  we will establish the inequality $\tilde b_{U,q} \leq b_{U,q}$  (in \cite {Ka2}, $U$ is a bounded symmetric domain  and $q=1$). 

For $r>0$, we set $U_r:=\{x\in U|\, |P(x)|\leq r\}$,  where  $P(x)$ is the {\it function}\/ on $U$ generated by the diastasis of the {\it original}\/ Kahler metric on $U$; $U_r$ is a compact  subset in $U$ \cite{T2, p.\,2, (1.2)}. Let $b(r)$ denote the $q$-Bergman metric on $U_r$. 

 It is well known that  $b_{U,q} = \lim_{r\to\infty} b(r)$. Given $U_r$, the restriction on $U_r$ of the projection
$U\rightarrow X_i$ is one-to-one for $i\gg 0$ because $\pi_1(X)$ is residually finite.  Further,  $b(r) > b_i|U_r$ for all $i> i(r)$. This  establishes the inequality $\tilde b_{U,q} \leq b_{U,q}$.

We have $\tilde b_{U,q}= \lim   b_i $. In the sequel, by a {\it local}\/ embedding we mean an embedding of a neighborhood. The above argument suggests the following proof.
\smallskip

(2.2) Let V be the completion of the Euclidean space $E:=\cup V_i$. We will show that $U$, with the metric $\tilde b_{U,q}$, can be naturally isometrically embedded in the infinite-dimensional projective space $\bold P(V^*)$ with its Fubini-Study metric. 

First, let $p(w)\in U$ be a point, and $ \Omega=f(dw_1\wedge \cdot\cdot\cdot \wedge dw_n)^q \in H$  in local coordinates $w=(w_1, \dots, w_n)$ around $p(w)$. Recall that the $b_{U,q}$-isometric embedding $U \hookrightarrow \bold P(H^*)$ was given by the $q$-Bergman kernel section $B:=B_{U, \Cal K^q}$ (see \cite {T1}; we do not assume $\Cal K_U^q$ is a trivial bundle). The map
$p(w)\mapsto \langle f,  B_w\rangle$ has defined
 the local $b_{U,q}$-isometric embedding $U \rightarrow \bold P(H^*)$ (see \cite {Kob}, \cite{FK, pp.\,5-13 or p.\,188}) hence the global $b_{U,q}$-isometric embedding (see \cite{C, Theorems 10, 11}).
\smallskip

(2.3) Now, let $\Omega =f(dw_1\wedge \cdot\cdot\cdot \wedge dw_n)^q$ be a form on $U$ that is the pullback of a form from $V_i$. As above, one can define locally the map $p(w)\mapsto \langle f,  B_w\rangle$ because $B =\lim_{r\to\infty} B(r)$ \cite{A, Part I, Sect.\,9}. 
The above map yields an analytic, local $b_i$-isometric embedding $U \rightarrow \bold P(V_i^*)$. Hence we get an analytic, local $\tilde b_{U,q}$-isometric embedding $U \rightarrow \bold P(V^*)$ (see \cite{FK, Proposition I.1.6 on p.\,13}). 

As in (2.2), by Calabi \cite{C, Theorems 10, 11}, the local $\tilde b_{U,q}$-isometric embedding yields the global one.
\smallskip

(2.4) Clearly, $U$ with the metric $b_{U,q}$ is $\bold P(H^*)$-resolvable at a point p (1-resolvable of rank $N=\infty$ in Calabi's terminology \cite{C, p.\, 19, Definition}); the image of $U$ does not lie in a proper subspace of $\bold P(H^*)$. 

Similarly, $U$ with the metric $\tilde b_{U,q}$ is $\bold P(V^*)$-resolvable at p. By the construction of the embedding $U\hookrightarrow \bold P(V^*) $, we get the embedding $\bold P(V^*) \hookrightarrow \bold P(H^*)$. It follows the latter embedding is surjective and $\bold P(V^*) = \bold P(H^*)$.

This concludes the proof of the proposition.
\enddemo

\proclaim{Corollary}  With the notation and assumptions of {\rm {(1.1)}},  $U$ is a
bounded domain in $\bold C^n$ provided $U$ has the $q$-Bergman metric for a sufficiently large integer $q$. 
\endproclaim

\demo{Proof} With all the assumptions of the corollary,  in the uniformization  theorem \cite{T1, Theorem} we can replace \lq\lq very large\rq\rq\/ by a {\it weaker}\/ assumption \lq\lq large\rq\rq. Indeed,  we now get the embedding corresponding to the metric $\tilde b_{U,q}$, as in \cite{T1, Sect.\,2.2}:
$$
 U \hookrightarrow {\bold P}(V^*).
$$

The above embedding is {\it required}\/  for the proof of the
uniformization theorem (see \cite{T1, Sect.\,3}); recall that the assumption \lq\lq {\it very}\/ large \rq\rq\/ (or \lq\lq $U$ has the $q$-Bergman metric\rq\rq)\/ does not follow from other assumptions \cite{T1, Remark 4.3}. 
\enddemo

\remark {Remark  \rm{1}} Let dim$X=1$ and $\Cal K_X$ is very ample. In the above proposition, one can take $m=1$, an arbitrary integer $q\geq m$, and  the Poincar\'e metric on the disk in place of the original metric. 

Other $\pi_1$-invariant volume forms on $U$ can be used in place of $dv_U$.
\endremark
\remark {Remark  \rm{2}} We {\it conjecture}\/ that, in the above corollary, one can replace   \lq\lq $U$ has the $q$-Bergman metric\rq\rq\/ by  \lq\lq the fundamental group $\pi_1(X)$ is nonamenable\rq\rq. 

Assume the fundamental group is nonamenable.  Perhaps, one has to employ the Bergman type reproducing kernel of a space of {\it harmonic}\/ forms. 
\endremark


\Refs
\widestnumber\key{ABC}

\ref \key A  \by N. Aronszajn   \pages 337--404
\paper Theory of Reproducing Kernels
\yr 1950   \vol 68
\jour Trans. Amer. Math .Soc.
\endref

\ref  \key  C  \by  E. Calabi \pages  1--23
\paper Isometric imbedding of complex manifolds
\yr1953 \vol  58
\jour Ann. of Math.
\endref

\ref  \key FK  \by  J. Faraut, S. Kaneyuki, A. Kor\'anyi, Q.-k. Lu, G. Roos
\book Analysis and geometry on complex
homogeneous domains
\publ Birkh\"auser, Boston
\yr 2000
\endref

\ref \key Ka1 \by D. A. Kazhdan \paper\nofrills On arithmetic varieties. 
\inbook  in Lie Groups and Their Representations  (Proc.  Summer School, Bolyai J\'anos Math. Soc., Budapest, 1971) 
\publaddr Halsted    \yr 1975   \pages  151--217
\endref

\ref \key Ka2 \bysame
\paper email to Curtis T.\, McMullen
\endref

\ref \key Kob  \by S. Kobayashi   \pages 267--290
\paper Geometry of bounded domains
\yr 1959   \vol 92
\jour Trans. Amer. Math .Soc.
\endref


\ref
\key Kol \by J. Koll\'ar
\book  Shafarevich maps and automorphic forms
\publ Princeton Univ. Press, Princeton
\yr 1995
\endref

\ref
\key M  \by  C.T. McMullen \pages
\paper Entropy on Riemann surfaces and the Jacobians of finite covers
\yr 
\toappear
\jour Comment. Math. Helv. 
\endref

\ref
\key O  \by  T. Ohsawa \pages
\paper Review and Questions on the Bergman Kernel in Complex Geometry
\yr 2010
\jour Preprint.
\endref

\ref
\key T1 \by  R. Treger \pages
\paper Uniformization
\yr
\jour arXiv:math.AG/1001.1951v4
 \endref

\ref
\key T2 \bysame  \pages
\paper Remark on a conjecture of Shafarevich
\yr
\jour arXiv:math.AG/1008.1745v2
 \endref

\ref
\key U \by  M. Umehara \pages 203-214
\paper Kaehler Submanifolds of Complex Space Forms
\yr 1987 \vol 10
 \endref

\ref
\key Y \by  S.-T. Yau \pages 108-159
\paper Nonlinear Analysis in Geometry
\yr 1987 \vol 33
\jour Eiseignement Math\'ematique
 \endref

\endRefs
\enddocument